\documentclass[preprint,3p,a4paper]{elsarticle}
\usepackage{amsmath,amssymb,amsthm}

\usepackage{hyperref}
\usepackage{geometry}
 \theoremstyle{plain}

  \newtheorem{thm}{Theorem}[section]
  \newtheorem{cor}[thm]{Corollary}
  \newtheorem{lem}[thm]{Lemma}
  \newtheorem{prop}[thm]{Proposition}
  \newtheorem{exmp}[thm]{Example}
  \newtheorem{rem}[thm]{Remark}
  
  \newtheorem{defn}[thm]{Definition}
  \newtheorem{prob}[thm]{Problem}

\newcommand{\da}{\downarrow \hspace{-2pt}}

\newcommand{\ua}{\uparrow \hspace{-2pt}}

\newcommand{\UUa}[0]{\makebox[1ex][l]{\lower.15ex
                                 \hbox{$\uparrow$}}\kern-1ex\lower-.15ex
                                 \hbox{$\uparrow$}}

\newcommand{\DDa}[0]{\makebox[1ex][l]{\lower.15ex
                                 \hbox{$\downarrow$}}\kern-1ex\lower-.15ex
                                 \hbox{$\downarrow$}}

\journal{Topology and its applications}
\begin{document}
\begin{frontmatter}
\title{\bf Continuity and core compactness of topological spaces\tnoteref{t1}}
\tnotetext[t1]{Research supported by NSF of China (Nos. 11871353, 12001385).}

\author{Yuxu Chen}
\ead{yuxuchen@stu.scu.edu.cn}
\author{Hui Kou\corref{cor}}
\ead{kouhui@scu.edu.cn}
\author{Zhenchao Lyu\corref{}}
\ead{zhenchaolyu@scu.edu.cn}
\cortext[cor]{Corresponding author}
\address{Department of Mathematics, Sichuan University, Chengdu 610064, China}

\begin{abstract}

We investigate two approximation relations on a $T_{0}$ topological space, the $n$-approximation and the $d$-approximation.\ Different kinds of continuous spaces are defined by the two approximations and are all shown to be directed spaces.\
We will show that the continuity of a directed space is very similar to the continuity of a dcpo in many aspects, which indicates that the notion of directed spaces is a suitable topological extension of dcpos. Besides, a series of characterizations of the core compactness of directed spaces are obtained, which are close related to a long-standing open problem that for which distributive continuous lattice its spectrum is exactly a sober locally compact Scott space.

  \vskip 3mm
{\bf Keywords}: approximation relation, continuity,  c-space, core compactness.

  \vskip 2mm

{\bf Mathematics Subject Classification}:  54A10, 54A20, 06B35.
\end{abstract}

\end{frontmatter}

\section{Introduction}

In domain theory, the continuity of a dcpo is defined through the waybelow relation. A dcpo is continuous iff the lattice of its Scott open subsets is a completely distributive lattice. More generally, any $T_{0}$ topological space is a c-space iff the lattice of its open subsets is a completely distributive lattice \cite{Erne91}. A continuous dcpo endowed with the Scott topology can be viewed as a special c-space. The natural question arises whether it is possible to define a approximation relation on $T_{0}$ topological spaces such that  c-spaces are some kind of continuous spaces defined through the approximation, similar to the continuity and the waybelow relation on a dcpo. 

Ershov \cite{ERS73} came up with an approximation relation on a topological space defined by that $x$ approximates $y $ iff $y \in (\ua x)^{\circ}$ in 1973. However, it is defined through the interior instead of directed subsets like the waybelow relation. In the past several years, some attempts were made to define approximation relations on a $T_{0}$ space like the waybelow relation on a dcpo. In 2014, Wang and Kou \cite{WK2014} introduced the $n$-approximation relation on $T_{0}$ topological spaces through convergent nets and then defined $n$-continuous spaces. They showed that a $T_0$ topological space is $n$-continuous iff it is a c-space. Thus, the $n$-approximation and $n$-continuity seem appropriate candidate for defining the continuity of a $T_{0}$ topological space.  However, these techniques in domain theory can not be used in $n$-continuous spaces since the $n$-approximation is based on nets rather than directed subsets.

Directed spaces \cite{YK2015} are special $T_{0}$ topological spaces,\ which are generalizations of Scott spaces. The idea is that the priori is a $T_0$ space and its convergent directed subsets relative to the specialization order rather than a poset and its existing directed supremums. In $T_{0}$ topological spaces, the notion of directed spaces is equivalent to that of monotone determined spaces defined by Ern\'{e}  \cite{ERNE2009}. $\bf DTop$, the category of directed spaces was shown to be cartesian closed.

$T_0$ c-spaces (or $n$-continuous spaces) are all directed spaces. Based on this, Feng and Kou \cite{FK2017} introduced another approximation relation on directed spaces, called $d$-approximation, and defined the quasicontinuous directed spaces. They showed that quasicontinuous spaces are same to locally hypercompact spaces. Recently, Zhang, Shi and Li \cite{ZHA2021} also came up with the notion of continuous spaces and quasicontinuous spaces. They considered these spaces restricted in weak monotone convergence spaces.

We will first investigate the relation between the two approximation relations and the continuity defined by them. A $T_{0}$ topological space is $n$-continuous (resp., $n$-algebraic, $n$-quasicontinuous, $n$-quasialgebraic) iff it is a $d$-continuous (resp., $d$-algebraic, $d$-quasicontinuous, $d$-quasialgebraic) directed space iff it is a c-space (resp., $b$-space, locally hypercompact space, hypercompactly based space). For convenience, we call them continuous (resp., algebraic, quasicontinuous, quasialgebraic) space uniformly. In classical domain theory, it is well known that continuous dcpos are the retracts of algebraic dcpos. It is the same for continuous spaces and algebraic spaces. Here, we give a proof by the structures called topological ideals, which are extensions of the rounded ideals of dcpos. 

When $X$ is a dcpo endowed with the Scott topology, it is a special directed space and the concepts of $d$-approximation and $d$-continuity agree with the waybelow relation and continuity in classical domain theory. The categorical products and exponential objects in the category of dcpos agree with those in $\bf DTop$ by viewing a dcpo as a directed space whose topology is the Scott topology. Thus, the notion of continuous spaces and directed spaces can be viewed as a natural extension of continuous domains and dcpos to $T_{0}$ topological spaces (see \cite{CK2020,FK2017,ZHA2021}). 

In domain theory, there is a correspondence between the continuity of a dcpo $L$ and the continuity of the lattice of its Scott open subsets $\sigma(L)$. More surprisingly, the continuity of a dcpo $L$ can be characterized by the continuity of its lattice of its Scott closed subsets $\Gamma(L)$ \cite{CK2022}. These correspondences can be extended to general $T_0$ spaces. We show that a $T_0$ space $X$ is continuous (resp., algebraic, quasicontinuous, quasialgebraic) iff the lattice of its closed subsets $C(X)$ is continuous (resp., algebraic, quasicontinuous, quasialgebraic).

A topological space is called core compact if its topology is a continuous lattice. There is a deep relation between core compactness and fuction spaces of topological spaces \cite{DK1970,ELS2004,Isbell1975}.  For example, a famous characterization is that   a $T_0$ space is core compactness  if and only if it is exponential in the category of $T_0$ spaces \cite
{CON03,GAU}. A poset $P$ endowed with the Scott topology $\sigma(P)$ is called a Scott space,  denoted by $\Sigma P$. When all topological spaces are restricted to the Scott spaces, core compactness can be characterized  simply by productions as follows: for a poset $P$, $\sigma(P)$ is a continuous lattice if and only if  for any poset $Q$ one has $\Sigma(P\times Q)=\Sigma(P)\times\Sigma(Q)$ \cite
{CON03}, which is equivalent to that  $\sigma(\Gamma(P)) $ is a continuous lattice with $\sigma(\Gamma(P))=\upsilon(\Gamma(P)) $ \cite{CK2022,MLZ2021}, where $\Gamma(P)$ is the lattice of all Scott closed subsets of $P$ and $\upsilon(\Gamma(P))$ is the upper topology of $\Gamma(P)$. These results show that the continuous Scott topologies of posets are a class of special distributively continuous lattices. Is there another class of $T_0$ space such that the continuous topology has the same  features? In this paper, we will further investigate core compactness of posets endowed with special topologies. Particularly, the results for core compact Scott spaces can be smoothly extended to directed spaces,  i.e., $T_0$ monotone determined spaces. These results are related to  a long-standing open problem that for which distributive continuous lattice its   spectrum is exactly a sober locally compact Scott space (see \cite[Problem 528]{Mill1990}).

\section{Preliminaries}

We assume some basic knowledge of domain theory and topology, as in, e.g., \cite{AJ94,CON03}.

A nonempty set $P$ endowed with a partial order $\leq$  is called a poset. For $A\subseteq P$, we set $\da A=\{x\in P: \exists a\in A,  \ x\leq a\}$,\ $\ua A =\{x\in P: \exists a\in P, \ a\leq x\}$.\ $A$ is called a lower or an upper set, if $A=\da A$ or $A=\ua A$ respectively. For an element $a\in P$, we use $\da a$ or $\ua a$ instead of $\da\{a\}$ or $\ua \{a\}$, respectively.

Let $P$ be a poset. We denote $\sigma(P)$, $\upsilon(P)$ and $A(P)$ to be the Scott topology,\ the upper topology and the Alexandroff topology on $P$ respectively. 
For a complete lattice $L$, and $x,y \in L$, we say that $x \prec y$ iff $y \in \mathrm{int}_{\upsilon(L)} \ua x$. $L$ is hypercontinuous if $\{d \in L:d \prec x\}$ is directed and has $x$ as a supremum for all $x\in L$. $L$ is hyperalgebraic if $x = \bigvee \{d\in L, d \prec d \leq x\}$ for all $x \in L$. An element $y$ such that $y \prec y$ is called hypercompact.

Topological spaces will always be supposed to be $T_0$. For a topological space $X$, its topology is denoted by $\mathcal{O}(X)$ or $\tau$. The partial order $\sqsubseteq$ defined on $X$ by $x\sqsubseteq y \Leftrightarrow  x\in \overline{\{y\}}$
is called the  specialization order, where $\overline{\{y\}}$ is the closure of $\{y\}$.\ From now on, all order-theoretical statements about  $T_0$ spaces, such as upper sets, lower sets, directed sets, and so on,  always refer to  the specialization order ``$\sqsubseteq $".\ 

For any two topological spaces $X,Y$, we denote $Y^{X}$ or $TOP(X,Y)$ the set of all continuous maps from $X$ to $Y$, endowed with the pointwise order. Denote $[X \to Y]_p$ and $[X \to Y]_I$ to be the topological space equipped with the topology of pointwise convergence and the Isbell topology on $Y^{X}$ respectively.  

\vskip 2mm

We now introduce the notion of a directed space. 

Let $(X,\mathcal{O}(X))$ be a $T_0$ space. Every directed subset $D\subseteq X$ can be  regarded as a monotone net $(d)_{d\in D}$. Set
$DS(X)=\{D\subseteq X: D \ {\rm is \ directed}\}$ to be the family of all directed subsets of $X$. For an $x\in X$, we denote  $D\rightarrow x$  to mean  that $x$ is a limit of $D$, i.e., $D$ converges to $x$ with respect to the topology on $X$. Then the following result is obvious.

\begin{lem}\rm  Let $X$ be a $T_0$ space. For any $(D,x)\in DS(X)\times X$, $D\rightarrow x$ if and only if $D\cap U\not=\emptyset$ for any open neighborhood of $x$.
\end{lem}

Set $DLim(X)=\{(D,x)\in DS(X)\times X: \ D\rightarrow x \}$ to be the set of all pairs of directed subsets and their limits in  $X$. Then $(\{y\},x)\in DLim(X) $ iff $x\sqsubseteq y$ for all $x,y\in X$. 

\begin{defn}\rm Let $X$ be a $T_0$ space. A subset $U\subseteq X$ is called directed-open  if for all $(D,x)\in DLim(X)$, $x\in U$ implies $D\cap U\not=\emptyset$.
\end{defn}

Obviously, every open set of $X$ is directed-open. Set
$d(\mathcal{O}(X))=\{U\subseteq X: U \ {\rm is} \  {\rm directed}\text{-}{\rm open}\},$
then $\mathcal{O}(X)\subseteq d(\mathcal{O}(X))$.

\begin{thm}\rm \cite{YK2015} \label{convergence}
 Let $X$ be a $T_0$ topological space. Then
\begin{enumerate}
\item[(1)] For all $U\in d(\mathcal{O}(X))$, $U=\ua U$;
\item[(2)] $X$ equipped with $d(\mathcal{O}(X))$ is a $T_0$ topological space such that  $\sqsubseteq_d =\sqsubseteq$, where $\sqsubseteq_d$ is the specialization order relative to $d(\mathcal{O}(X))$.
\item[(3)] For a directed subset $D$ of $X$, $D\rightarrow x$ iff $D\rightarrow_d x$ for all $x\in X$, where $D\rightarrow_d x$ means that $D$ converges to $x$ with respect to the topology $d(\mathcal{O}(X))$.

\item[(4)] $d(d(\mathcal{O}(X))=d(\mathcal{O}(X))$.
\end{enumerate}
\end{thm}

\begin{defn}\rm \cite{YK2015} A  topological space $X$ is said to be a directed space if it is $T_0$ and  every directed-open set is open; equivalently, $d(\mathcal{O}(X))=\mathcal{O}(X)$.
\end{defn}

One can see that the idea to define a directed space is similar to define a sequential space and the Scott topology on a poset. In $T_{0}$ topological spaces, the notion of a directed space is equivalent to a monotone determined space defined by Ern\'e \cite{ERNE2009}.\  Given any space $X$, we denote $\mathcal{D}X$ to be the topological space $(X,d(\mathcal{O}(X)))$.

\begin{thm}\rm \cite{YK2015} Let $X$ be a $T_0$ space, we have
\begin{enumerate}
\item[(1)] $\mathcal{D}X$ is a directed space.

\item[(2)] The following three conditions are equivalent to each other:
\begin{enumerate}
\item[(i)] $X$ is a directed space;
\item[(ii)] For all $U\subseteq X$, $U$ is open iff for any $(D,x)\in DLim(X)$, $x\in U$ implies $U\cap D\not=\emptyset$.
\item[(iii)] For all $A\subseteq X$, $A$ is closed iff for any directed subset $D\subseteq A$, $D\rightarrow x$ implies $x\in A$ for all $x\in X$.
\end{enumerate}
\end{enumerate}
\end{thm}

Directed spaces include many important structures in domain theory.

\begin{exmp}\rm \

\begin{enumerate}[(1)]
\item Every poset endowed with the Scott topology is a directed space.

\item (\cite{ERNE2009})Every poset endowed with the weak Scott topology is a directed space.
\item Every c-space is a directed space. In particular, any poset endowed with the Alexandroff topology is a directed space.
\item (\cite{ERNE2009,FK2017}) Every locally hypercompact space is a directed space.
\end{enumerate}

\end{exmp}

We denote ${\bf DTop}$ to be the category of all nonempty directed spaces with continuous maps as morphisms. Then the following statements hold.

\begin{thm}\rm \cite{L2017,YK2015} $\bf DTop $ is cartesian closed. Let $X,Y$ be directed spaces. 

\begin{enumerate}[(1)]
\item The categorical product $X \otimes Y$ of $X$ and $Y$ is homeomorphic to $\mathcal{D} (X \times Y)$.

\item The exponential object $[X \to Y]$  is homeomorphic to $ \mathcal{D}([X \to Y]_p)$.

\end{enumerate}

\end{thm}

\vskip 3mm

\begin{lem} \rm  \label{directed product}
Given any two directed spaces $X,Y$, a subset $U$ of $X \otimes Y$ is open iff the following conditions hold:
\begin{enumerate}[(1)]
\item for any directed subsets $(x_{i})_{I}$ with $(x_{i})_I \to x$ in $X$ and any $(x,y) \in U$, we have $((x_{i},y))_{I} \cap U \neq \emptyset$;

\item for any directed subsets  $(y_{i})_{I}$ with $(y_{i})_I \to y$ in $Y$ and any $(x,y) \in U$ we have $((x,y_{i}))_{I} \cap U \neq \emptyset$.
\end{enumerate}
\end{lem}

\noindent{\bf Proof.} Assume that $U$ is open in $X \otimes Y$ and $(x,y) \in U$. For any directed set $(x_{i})_I \to x$ in $X$ and any $y \in Y$, $((x_{i},y))_{I}$ is a directed subset of $X \otimes Y$ and $((x_{i},y))_{I} \to (x,y)$ in $X \times Y$. By Theorem \ref{convergence} (3), $((x_{i},y))_{I} \to (x,y)$ in $X \otimes Y$. Thus, $((x_{i},y))_I \cap U \neq \emptyset$. It is the same for $(y_{i})_{I}$.

Conversely, assume that $(1)$ and $(2)$ are satisfied. We show that $U$ is open in $X \otimes Y$. It is easily seen that $U$ is an upper set relative to the specialization order of $X \times Y$. 
Let $D = ((x_{i},y_{i}))_{I}$ be a directed subset of $X \times Y$ and converges to $(x,y) \in U$ in $X \times Y$. We have that $(x_{i})_I \to x$ in $X$ and $(y_i)_I \to y$ in $Y$ respectively. Thus, $((x_{i},y))_{I} \to (x,y)$ in $X \times Y$, and then there exists some $i_{0} \in I$ such that $(x_{i_{0}},y) \in U $.\ By $((x_{i_{0}}, y_{i}))_{I} \to (x_{i_{0}},y)$, there exists some $i_{1} \in I$ such that $(x_{i_{0}},y_{i_{1}}) \in U$. Let $i_{0},i_{1} \leq i_{2}$, then $(x_{i_{2}},y_{i_{2}}) \in U$. $\Box$

\section{Approximation relations and continuity of topological spaces}

In this section, we will first investigate the relationships between the $n$-approximation and the $d$-approximation on a $T_{0}$ topological space.\ We will see that  a $T_{0}$ topological space is $n$-continuous (resp., $n$-algebraic) iff it is a $d$-continuous (resp., $d$-algebraic) directed space iff it is a c-space (resp., $b$-space).  These results indicate that the continuity defined by $n$-approximation and $d$-approximation is a suitable extension of continuity of a dcpo to topological spaces. Continuous spaces in directed spaces are just as continuous domains in dcpos.

\vskip 3mm

We first compare the $n$-approximation and the $d$-approximation.

\begin{defn} \rm \cite{FK2017,WK2014}
Let $X$ be a $T_{0}$ topological space and $x,y\in X$. 
\begin{enumerate}[(1)]
\item We say that $x$ $d$-approximates $y$, denoted by $x\ll_d y$, if for any directed subset $D\subseteq X$, $D\rightarrow y$ implies $x\sqsubseteq d$ for some $d\in D$. If $x\ll_d x$, then $x$ is called a $d$-compact element of $X$.

\item We say that $x$  $n$-approximates $y$ if for all nets $(y_j)_{j\in J}$ (shortly, $(y_j)_{J}$) of $X$, $(y_j)_{J}\rightarrow y$ implies $x\sqsubseteq y_{j_0}$ for some $j_0\in J$. We call $x$ $n$-compact if it $n$-approximates itself.
\end{enumerate}

\end{defn}

We introduce the following notations for $x,y\in X$, where $i = d,n$:
\begin{eqnarray*}
\DDa_i x &=& \{y\in X: y\ll_i x\}\\
\UUa_i x &=& \{y\in X: x\ll_i y\}\\
K_i(X)  &=& \{x\in X: x\ll_i x\}
\end{eqnarray*}

It is easily seen that the definition of $n$-approximation is stronger than that of $d$-approximation, i.e.,  $x \ll_{n} y$ implies $x \ll_{d} y$. Is the converse true? We give two negative examples, where the first one is endowed with the upper topology and the second one is endowed with the Scott topology and hence it is a directed space.

\begin{exmp} \rm  \label{exand}
Let $\mathbb{N}$ be the set of natural numbers.
\begin{enumerate}[(1)]

\item  Denote $\mathbb{N}^\top$ the flat domain, i.e., the poset with carrier set $\mathbb{N} \cup \{\top\}$ and $x \leq y$ iff $y = \top$ or $x=y$. 
It is easily seen that $\mathbb{N}^\top$ is a dcpo. 
Denote $\mathbb{N}^\top_\bot$ be the $\mathbb{N}^\top$ adding a bottom element, then it is a complete lattice.

\item Let $\omega = \{\omega_{i}: i \in \mathbb{N}\} $.  Let $ P = (\mathbb{N} \times \mathbb{N}) \cup \ \omega $ endowed with a partial order as follows:

\begin{enumerate}

\item[(i)] $\forall (m_{1},n_{1}),(m_2,n_2) \in \mathbb{N} \times \mathbb{N} $, $(m_1,n_1) \leq (m_2,n_2)$ iff $m_1 = m_2 $ and $ n_{1} \leq n_{2}$;
\item[(ii)] $\forall i,m,n \in \mathbb{N}$, $(m,n) \leq \omega_{i}$ iff $m =i$ or $m=1, n \leq i$.

\end{enumerate}
\end{enumerate}

\begin{enumerate}
  \item[(1)]  Considering the upper topology on $\mathbb{N}^\top$ (or $\mathbb{N}^\top_\bot$), we have: 
\begin{enumerate}
\item $\top \ll_d \top$; $\top \not \ll_n \top$;
\item $(\mathbb{N}^\top, v(\mathbb{N}^\top))$ is a locally compact sober space;
\item  $(\mathbb{N}^\top, v(\mathbb{N}^\top))$ is not a directed space.
\end{enumerate}

 \item[(2)]  Considering $P$ endowed with the Scott topology, we have $\forall n \in \mathbb{N},\ (1,n) \ll_d \omega_1$ \& $ (1,n) \not\ll_n \omega_1$.

\end{enumerate}

\end{exmp}

From the two examples we know that even in directed spaces, the two approximation relations may be different.

\begin{lem} \rm \label{topo ideal}
Let $X$ be a $T_{0}$ topological space and $D$ be a directed subset of $X$.\
\begin{enumerate}[(1)]
\item If $D\subseteq\da x$ and $D\rightarrow x$ for $x\in P$, then $x=\sup D$.
\item If $D \to x$ and $\sup D$ exists,\ then $x \leq \sup D$.

\end{enumerate}
\end{lem}

\noindent{\bf Proof.}
(1) Suppose that $D\subseteq\da x$ and $D\rightarrow x$ for $x\in X$.\ Assume that there exists $y\in X$ such that $D\subseteq \da y$ and $x\not\leq y$.\ Then $x\in P\backslash\!\da y\in \mathcal{O}(X)$. Thus $D\cap (X\backslash\!\da y ) \not=\emptyset$.\ It is a contradiction.\ Hence,\ $x=\sup D$.\

(2) Since $D \to x$,\ given any open subset $U$ such that $x \in U$,\ there exists some $d \in D$ such that $d \in U$ and then $\sup D \in U$.\ Therefore,\ $x \leq \sup D$. $\Box$

\vskip 2mm

By Lemma \ref{topo ideal}, we know that given any dcpo $P$ and any directed subset $D$ of $P$, then $D \to x \in P$ in $\sigma(P)$ iff $x \leq \sup D$. Thus, when $P$ is a dcpo, the $d$-approximation on $(P,\sigma(P))$ is equal to the waybwlow relation on $P$. By Example \ref{exand} (2), the $n$-approximation may be different to the waybelow relation in a dcpo. The notion of $d$-approximation seems better than $n$-approximation as an extension of waybelow relation. 
We now introduce the definition of continuous spaces and we will see that for a dcpo $P$, $(P,\sigma(P))$ is a continuous space if and only if $P$ is a continuous dcpo in the sense of the classical domain theory. Moreover, in any continuous space, the $d$-approximation is equal to $n$-approximation.

\begin{defn} \rm
Let $X$ be a $T_{0}$ topological space.

\begin{enumerate}[(1)]
\item $X$ is called $d$-$continuous$ if it is a directed space such that $\DDa_dx$ is directed and $\DDa_dx\rightarrow x$ for all $x\in X$.

\item $X$ is called $n$-$continuous$ if for every $x\in X$, there exists a net $(x_j)_J\subseteq \DDa_{n} x$ such that $(x_j)\rightarrow x$.
\end{enumerate}

\end{defn}

\begin{rem} \rm
When defining a $d$-continuous space, we require it to be a directed space. Otherwise, it may not be a c-space. Example \ref{exand} (1) is an example. $n$-continuous spaces and c-spaces are all directed spaces.
Therefore, it would better restrict the spaces in directed spaces when discussing continuity.
\end{rem}

\begin{lem}\rm \cite{CK2020} \label{d-continuity}
Let $X$ be a $d$-continuous space. Then we have the following statements.
\begin{enumerate}[(1)]
\item For all $x,y\in X$, $x\ll_d y$ implies $x\ll_d z\ll_d y$ for some $z\in X$.
\item $\{\UUa_d x: x\in X\}$ is a basis of the topology of $X$.

\item For all $x,y\in X$, the following are equivalent:
\begin{enumerate}
\item[(i)] $x\ll_d y$;
\item[(ii)] $y\in (\ua x)^{\circ}$.
\item[(iii)] For any net $(x_j)_J\subseteq X$, $(x_j)\rightarrow y$ implies $x\sqsubseteq x_{j_0}$ for some $j_0\in J$.

\end{enumerate}
\end{enumerate}

\end{lem}

\begin{lem}\rm \cite{WK2014}\label{n-continuity}
For a $T_0$ space $X$, the following are equivalent: for all $x,y\in Y$,
\begin{enumerate}
\item[(1)] $x\ll_n y$;
\item[(2)] $y\in (\ua x)^{\circ}$, where  $(\ua x)^{\circ}$ means the interior of $\ua x$ ;
\item[(3)] For all net $(y_j)_J$, $(y_j)\rightarrow y$ implies $\{y_j\}_J$ is eventually in $\ua x$.

\end{enumerate}
\end{lem}

We can see from Lemma $\ref{n-continuity}$ that $x \ll_{n} y$ iff $x$ approximates $y$ in the sense of Ershov.\ Thus,\ the notion of $n$-continuous spaces is in fact equivalent to that of c-spaces. For the concepts of continuity in domain theory, the notion of a c-space is another famous generalization.

\begin{defn}\rm
A topological space $X$ is called a c-space, if it is $T_0$ and every point $y$ has a neighborhood basis of sets of the form $\ua x$.
\end{defn}

The c-spaces had been studied by Ern\'{e} \cite{Erne81,Erne91}, and by Ershov \cite{Ershov93,Ershov97} under the name of an $\alpha$-space. 
Adding the notions of $d$-continuity and $n$-continuity we introduce above, there are three generalizations of continuity in domain theory. Combining the results in \cite{CK2020} and \cite{WK2014}, we have the following statement.

\begin{thm}\rm \cite{CK2020,WK2014}  \label{continuous space}
For a $T_{0}$ topological space $X$, the following are equivalent:
\begin{enumerate}
\item[(1)] $X$ is a $d$-continuous directed space;
\item[(2)] $X$ is a directed space and for each $x\in X$, there is a directed subset $D\subseteq \DDa_d x$ such that $D\rightarrow x$;
\item[(3)] $X$ is  $n$-continuous;
\item[(4)] for every $x\in X$, $\DDa_n x$ is directed and converges to $x$;
\item[(5)] for every $x\in X$, there exists a  directed subset $D \subseteq \DDa_n x$ such that $ D \to x$;

\item[(6)] $X$ is a c-space;
\item[(7)] $\mathcal{O}(X)$ is a completely distributive lattice.

\end{enumerate}
\end{thm}
\noindent{\bf Proof.}
$(3) \Leftrightarrow (4) \Leftrightarrow (6)  $ was proved in \cite{WK2014}. $(6) \Leftrightarrow (1) \Leftrightarrow (2)$ was proved in \cite{CK2020}. $(4) \Rightarrow (5)$ and $(5) \Rightarrow (3)$ are obvious. $\Box$

\vskip 2mm

We give a uniform name to replace the three notions of continuity.

\begin{defn}\rm \label{definition of continuous space}
A topological space is called $continuous$ if one of the seven equivalent conditions in Theorem \ref{continuous space} holds.
\end{defn}

From Lemma \ref{d-continuity} and Lemma \ref{n-continuity}, the $d$-approximation and $n$-approximation coincide when $X$ is a continuous space. For convenience, both $d$-approximation and $n$-approximation  in a continuous space are replaced by $\ll$ uniformly in the sense of Definition \ref{definition of continuous space}. From now on, we do not use the subscripts $d,n$ when no confusion can arise.

\begin{rem}\rm
\begin{enumerate}
\item[(1)] For not leading to confusion, a continuous dcpo or poset  always means the order theoretical one in domain theory.

\item[(2)] Every continuous poset is exactly a continuous space whose topology is the Scott topology.

\item[(3)] Every $T_0$ Alexandroff space is a continuous space. So continuity in our sense is a proper extension of the one in domain theory.
\end{enumerate}
\end{rem}

We next give the definition of an algebraic space.

\begin{defn}\rm
Let $X$ be a topological space. $X$ is called $d$-algebraic if it is a directed space such that $\da x \cap K_{d}(X)$ is directed and converges to $x$ for all $x\in X$. $X$ is called $n$-algebraic if $\da x \cap K_{n}(X)$ is directed and converges to $x$ for all $x\in X$.
\end{defn}

Obviously, all $d$-algebraic spaces and $n$-algebraic spaces are continuous spaces. By Theorem \ref{d-continuity} and Theorem \ref{n-continuity}, we know that $d$-compact elements coincide with $n$-compact elements for continuous spaces. Thus, the notion of a $d$-algebraic space is equivalent to a $n$-algebraic space. We call them algebraic spaces uniformly. It is easy to see that a space is algebraic
iff it is a b-space. Since a space is a b-space iff the lattice of its open subsets is an completely distributive algebraic lattice \cite{Erne91}, we have the following characterizations for algebraic spaces. 

\begin{thm}\rm
For a $T_{0}$ topological space $X$, the following conditions are equivalent:
\begin{enumerate}
\item[(1)] $X$ is algebraic;
\item[(2)] $X$ is a b-space, i.e., $\forall x \in U \in \mathcal{O}(X), \exists y \in U.\ x \in (\ua y)^{\circ} = \ua y$;
\item[(3)] $\mathcal{O}(X)$ is a completely distributive algebraic lattice.
\end{enumerate}
\end{thm}

When $X$ is a poset with the Scott topology, the above notion of an algebraic space coincides with the order theoretical one in domain theory. In classical domain theory, it is shown that continuous domains are the retracts of algebraic domains \cite{CON03}. The result can be extended to continuous spaces as well. We will construct a canonical algebraic space for each $T_0$ space and show that a $T_0$ space is continuous if and only if it is a continuous retract of an algebraic space.

\begin{defn} \rm Let $X$ be a $T_0$ space. A directed subset $D\subseteq X$ is called an ideal net if there exists  $x\in X$ such that:
\begin{enumerate}
\item[(1)] $D\rightarrow x$,
\item[(2)] $\forall d\in D$, $d\sqsubseteq x$.
\end{enumerate}
\end{defn}

Obviously, if $X$ is a poset endowed with the Scott topology, then every directed subset with a existing supremum is an ideal net.

\begin{lem} \rm If $D\subseteq X$ is an ideal net of a  $T_0$ space  $X$, then $\bigvee D$ exists and $D\rightarrow \bigvee D$.
\end{lem}

\noindent{\bf Proof.} By Lemma \ref{topo ideal}.

\begin{defn} \rm Let $X$ be a $T_0$ space. A subset $A$ of $X$ is called a topological ideal if there is an ideal net $D\subseteq X$ such that $A=\da D$. \end{defn}

So a topological ideal $A$ is a directed lower subset such that it has a supremum (denoted by $\bigvee A$) and converges to its supremum. Particularly, we set ${\rm I}_T(X)$ to be the set of all topological ideals of a $T_0$ space $X$.  Obviously, we have  $\{\da x: x\in X\}\subseteq {\rm I}_T(X)$. Set
$$\Omega({\rm I}_T(X))=\{{\mathcal U}\subseteq {\rm I}_T(X):  A\in {\mathcal U}\Leftrightarrow  \exists x\in X, \ \da x\in {\mathcal U} \ \& \da x\subseteq A\}.$$
Particularly, for any $x\in X$, we set ${\mathcal U}_x=\{A\in {\rm I}_T(X): x\in A\}$, then  ${\mathcal U}_x\in \Omega({\rm I}_T(X))$.

\begin{lem} \rm 
Let $X$ be a $T_0$ space. Then $({\rm I}_T(X),\Omega({\rm I}_T(X)))$ is an algebraic space such that:
\begin{enumerate}
\item[(1)] Its specialization order is equal to the set-theoretical inclusion;
\item[(2)] $K({\rm I}_T(X))=\{\da x: x\in X\}$.
\end{enumerate}
\end{lem}

\noindent{\bf Proof.}  Obviously, $\emptyset, {\rm I}_T(X) \in\Omega({\rm I}_T(X))$ and $\Omega({\rm I}_T(X))$ is closed under unions.  For ${\mathcal U},{\mathcal V}\in \Omega({\rm I}_T(X))$ and  $A\in {\mathcal U}\cap {\mathcal V}$,  there are $x,y\in X$ such that $\da x\in \mathcal U, \ \da y\in \mathcal V$ and $x,y\in A$. Since $A$ is directed, there exists $a\in A$ with $x,y\sqsubseteq a$. Thus $\da a\in \mathcal U\cap \mathcal V$ and $\da a\subseteq A$, i.e., $\mathcal U\cap \mathcal V\in \Omega({\rm I}_T(X))$. Hence, $\Omega({\rm I}_T(X))$ is a topology. If  $A\not\subseteq B$ for $A,B\in {\rm I}_T(X)$, then there is $a\in A$ together with $a\not\in B$; thus $A\in {\mathcal U}_a$ and $B\not\in{\mathcal U}_a$. If  $A \subseteq B$ for $A,B\in {\rm I}_T(X)$ and $A \in \mathcal U \in \Omega({\rm I}_T(X))$, then there exists some $\da x$ such that $\da x \subseteq A , \da x \in \mathcal U$; we have $\da x \subseteq B$ and thus $B \in \mathcal{U}$. Hence, $({\rm I}_T(X),\Omega({\rm I}_T(X)))$ is  $T_0$ and its specialization order is equal to the set-theoretical inclusion. 

Let $x \in X$, $\mathcal B$ be a directed subset of ${\rm I}_T(X)$ and $\mathcal{B} \to \ \da x$. Since $\mathcal U_{x}$ is open, there must exists some $B \in \mathcal B $ such that $B \in \mathcal U_{x}$, i.e., $\da x \subseteq  B $. Hence, we have that $\da x$ is a compact element of $({\rm I}_T(X), \Omega({\rm I}_T(X))))$. For any $A\in {\rm I}_T(X)$, $\mathcal{A}=\{\da x:x\in A\}$ is directed. For $\mathcal{U}\in \Omega({\rm I}_T(X))$ with $A\in\mathcal{U}$, there exists $a\in A$ such that $\da a\in \mathcal{U}$. Thus,  $\mathcal{A}$ converges to $A$ in ${\rm I}_T(X)$. It follows that $({\rm I}_T(X),\Omega({\rm I}_T(X)))$ is an algebraic space and  $K({\rm I}_T(X))=\{\da x: x\in X\}$. $\Box$

\vskip 3mm

By definition, every topological ideal is directed and converges to its supremum. Hence, there exists a surjective map from ${\rm I}_T(X)$ onto $X$ as follows:
$$\bigvee: {\rm I}_T(X)\longrightarrow X, \ \ \forall A\in {\rm I}_T(X), \ A\mapsto \bigvee A.$$

\begin{prop} \rm
For a $T_0$ space $X$,  $\bigvee: {\rm I}_T(X)\longrightarrow X$ is a continuous surjective map.
\end{prop}

\noindent{\bf Proof.} Obviously, $\bigvee: {\rm I}_T(X)\longrightarrow X$ is surjective. For any $U\in \mathcal{O}(X)$ and for $A\in\bigvee^{-1}(U)$,
we have $ \bigvee A\in U$. Note that since $A$ is a topological ideal, we have that $A=\da A$ and $A\rightarrow \bigvee A$ in $X$. There exists $a\in A$ such that $a\in U$. Thus, $\da a\in \bigvee^{-1}(U)$ and $\da a\subseteq A$. Hence,  $\bigvee^{-1}(U)\in \Omega({\rm I}_T(X))) $, i.e., $\bigvee: {\rm I}_T(X)\longrightarrow X$ is continuous. $\Box$

\vskip 3mm

Hence, every $T_0$ space is a continuous image of some algebraic space. Particularly, when $X$ is a dcpo endowed with the Scott topology $\sigma(X)$,  ${\rm I}_T(X)$ is exactly the ideal completion of $X$.

\begin{prop} \rm
Let $X$ be a dcpo endowed with the Scott topology $\sigma(X)$. Then  ${\rm I}_T(X)=Id(X)$, where $Id(X)$ is the set of all  ideal of  $X$ endowed with the inclusion order.
\end{prop}
\noindent{\bf Proof.} Given any topological ideal $D$, it is an ideal by definition. Thus, $D \in Id(X)$. Conversely, suppose that $A$ is a ideal of $X$. Since $X$ is a dcpo endowed with the Scott topology and $A$ is a directed subset, we know that $A \to \bigvee A$. Hence, $A$ is a topological ideal. $\Box$

\vskip 2mm
The following is the definition of an adjunction between posets.

\begin{defn} \rm
\cite{CON03}
Let $f:P\longrightarrow Q$ and  $g: Q\longrightarrow P$ be a pair of monotone maps between posets $P$ and $Q$. We say that $(f,g)$ is an adjunction provided the following holds:
$$\forall x\in P,y\in Q, \ f(x)\leq y\Leftrightarrow x\leq g(y).$$
In this case, $f$ is called  the lower adjoint of  $g$ and  $g$ is called the upper adjoint of $f$.
\end{defn}

\begin{thm} \rm
Let $X$ be a  $T_0$ space. Then following conditions are equivalent to each other:
\begin{enumerate}
\item[(1)] $X$ is a continuous space;
\item[(2)] $\bigvee: {\rm I}_T(X)\longrightarrow X$ has a continuous lower adjoint;
\item[(3)] $X$ is a continuous retract of some algebraic space.
\end{enumerate}
\end{thm}

\noindent{\bf Proof.} $(1)\Rightarrow (2)$. Let  $X$ is a continuous space. Then $\DDa x$ is directed and converges to $x$ for any $x\in X$ by Theorem \ref{continuous space}. Note that since $\DDa x=\da(\DDa x)$, we have that  $\DDa x$ is a topological ideal. Define a mapping  $\DDa: X\longrightarrow {\rm I}_T(X)$ as follows:  $$\forall x\in X, \ \DDa (x)=\DDa x.$$
Then $\DDa: X\longrightarrow {\rm I}_T(X)$ is well defined and monotone. For any $a\in X$ and for  $x\in \DDa^{-1}({\mathcal U}_a)$, we have  $a\in \DDa x$. Thus, there is some $y\in X$ such that  $a\ll y\ll x$. Set  $U=\UUa y$. Then  $U$ is an open neighborhood of  $x$ with  $U\subseteq \DDa^{-1}({\mathcal U}_a)$. Hence, $\DDa^{-1}({\mathcal U}_a)$ is an open subset of  $X$. Note that since $\{{\mathcal U}_a:a\in X\}$ is a base of  $\Omega({\rm I}_T(X))$, we have that  $\DDa: X\longrightarrow {\rm I}_T(X)$ is a continuous map. For any $x\in X$ and $A\in {\rm I}_T(X)$, $x\leq \bigvee A$ if and only if  $\DDa x\subseteq A$. Therefore,  $\DDa: X\longrightarrow {\rm I}_T(X)$ is the lower adjoint of  $\bigvee: {\rm I}_T(X)\longrightarrow X$ by definition.

$(2)\Rightarrow (3)$. Let $f:X\longrightarrow {\rm I}_T(X)$ be a continuous lower adjoint of $\bigvee: {\rm I}_T(X)\longrightarrow X$. Then  $\bigvee\circ f=id_X$ for that $\bigvee$ is surjective. Hence, $X$ is a continuous retract of the algebraic space ${\rm I}_T(X)$.

$(3)\Rightarrow (1)$. Let $Y$ be an algebraic space such that $X$ is a continuous retract of $Y$. Then there exist two continuous maps  $f:Y\longrightarrow X$ and  $g:X\longrightarrow Y$ with  $f\circ g=id_X$. For any $x\in X$ and $k\in K(Y)$, if  $k\leq g(x)$, then $f(k)\leq f(g(x))=x$. Let  $(x_j)_{j\in J}$ be a net of  $X$ with $(x_j)_J\rightarrow x$. Since  $g$ is continuous, we have  $(g(x_j))_J\rightarrow g(x)$. As  $k$ is compact and  $k\leq g(x)$, there exists  $j_0\in J$ such that  $k\leq g(x_{j_0})$. Hence, $f(k)\leq f(g(x_{j_0}))=x_{j_0}$, i.e., $f(k)\ll x$ by definition. Note that since  $Y$ is an algebraic space, $\da g(x)\cap K(Y)$ is directed and converges to $g(x)$. Thus, $\{f(z):z\in \da g(x) \cap K(Y)\}\subseteq \DDa x$ and $(f(z))_{z \in \downarrow g(x) \cap K(Y)}\rightarrow f(g(x))=x$. Therefore, $X$ is a continuous space. $\Box$

\vskip 3mm

\section{Continuity of lattices of closed subsets}

Similarly with continuous spaces and algebraic spaces in the previous section, we introduce the concepts of quasicontinuous spaces and quasialgebraic spaces.\ $n$-quasicontinuous spaces are equal to $d$-quasicontinuous spaces as well and they are equal to locally hypercompact spaces. We will summarize the correspondences among the continuity of a $T_0$ space, the continuity of the lattice of its open subsets, and the continuity of the lattice of its closed subsets.

Given any set $X$, we use $F \subseteq_{f} X$ to denote that $F$ is a finite subset of $X$.\ Given any two subsets $G,H \subseteq X$,\ we define $G \leq H$ iff $\ua H \subseteq \ua G$.\ A family of finite sets is said to be directed if given $F_{1}, F_{2}$ in the family,\ there exists a $F$ in the family such that $F_{1},F_{2} \leq F$.\

\begin{defn} \rm \cite{FK2017}
Let $X$ be a directed space, $y \in X$ and $G,H \subseteq X$.\ We say that $G$ $d$-$approximates$ $H$, denoted by $G \ll_{d} H$, if for every directed subset $D\subseteq X$, $D \to y$ for some $y\in H$ implies $D\ \cap \ua\! G \neq \emptyset$. We write $G \ll_{d} x$ for $G \ll_{d} \{ x \}$. $G $ is said to be compact if $G \ll_{d} G $.
\end{defn}

Let $X$ be a $T_{0}$ topological space and $\mathcal{M}$ be a directed family of finite subsets of $X$.\ Then $F_{\mathcal{M}} = \{A \subseteq X : \exists M \in \mathcal{M}, \ua M \subseteq A\}$ is a filter under the inclusion order.\ We say that $\mathcal{M}$ converges to $x$ if $F_{\mathcal{M}}$ converges to $x$,\ i.e.,\ for any open neighbourhood $U$ of $x$,\ there exists some $A \in \mathcal{M}$ such that $A \subseteq U$.\

\begin{defn} \rm \cite{FK2017}
A topological space $X$ is called $d$-$quasicontinuous$ if it is a directed space such that for any $x \in X$,\ the family $fin_{d}(x) = \{F:\ F  $ is finite,\ $F \ll_{d} x \}$ is a directed family and converges to $x$.
\end{defn}

Obviously, every continuous space is $d$-quasicontinuous. When $P$ is a poset endowed with the Scott topology $\sigma(P)$, $P$ is $d$-quasicontinuous if and only if $P$ is quasicontinuous in the sense of the classical domain theory.\ Given any quasicontinuous space $X$ and $F \subseteq_{f} X$,\ we write $\Uparrow_{d} F = \{ x \in X:\ F \ll_{d} x \}$.\

\begin{lem} \rm \cite{FK2017}\label{d-quasicontinuity}
Let $X$ be a $d$-quasicontinuous space. Then we have the following.
\begin{enumerate}
\item[(1)] Given any $H \subseteq X$ and  $y\in X$, $H \ll_d y$ implies $H \ll_d F\ll_d y$ for some finite subset $F \subseteq_{f} X$.

\item[(2)] Given any $F \subseteq_{f} X $, $\Uparrow_d F =(\ua F)^{\circ}$. Moreover, $\{\Uparrow_d F: F \subseteq_{f} X\}$ is a basis of the topology of $X$.

\item[(3)] Given any $G,H \subseteq X$, the following are equivalent:
\begin{enumerate}
\item[(i)] $G\ll_d H$;
\item[(ii)] $H \subseteq (\ua G)^{\circ}$.
\item[(iii)] For any net $(x_j)_J\subseteq X$, $(x_j)\rightarrow y, y \in H $ implies $\ua G \cap  (x_{j})_{J} \neq \emptyset$.

\end{enumerate}
\end{enumerate}

\end{lem}

\noindent{\bf Proof.}
 $(1)$ and $(2)$ was proved in \cite{FK2017} and the proof is analogous to that for quasicontinuous domains. We only prove $(3)$.

(i)$\Rightarrow $(ii) follows directly from (1), (2) and the fact that $G \ll_{d} H$ iff $G \ll_{d} x$ for every $x \in H$. (ii)$\Rightarrow $(iii) and (iii)$\Rightarrow $(i) are obvious from the definitions of convergence and $d$-approximation.
$\Box$

\vskip 2mm

Now we define $n$-quasicontinuous spaces.
\begin{defn} \rm
Let $X$ be a $T_{0}$ space and $G,H \subseteq X$.\ We say that $G$ $n$-$approximates$ $H$, denoted by $G \ll_{n} H$, if for every net $(x_{j})_{J}$ of $X$, $(x_{j})_{J} \to y$ for some $y\in H$ implies $(x_{j})_{J} \cap \ua G \neq \emptyset$. We write $G \ll_{n} x$ for $G \ll_{n} \{ x \}$. $G $ is said to be $n$-compact if $G$ is a finite set and $G \ll_{n} G $.
\end{defn}

\begin{lem}\rm \label{n-quasicontinuity}
For any $T_0$ space $X$, the following are equivalent: for all $x,y\in Y$,
\begin{enumerate}
\item[(1)] $G\ll_n H$;
\item[(2)] $H \subseteq (\ua G)^{\circ}$, where  $(\ua G)^{\circ}$ means the interior of $\ua G$ ;
\item[(3)] For every net $(y_j)_J$, $(y_j)\rightarrow y$ and $y \in H$ implies that $\{y_j\}_J$ is eventually in $\ua G$.

\end{enumerate}
\end{lem}

{\noindent\bf Proof.} $(2)\Rightarrow (3)$ and $(3)\Rightarrow (1)$ are obvious. We only show $(1)\Rightarrow (2)$. Assume that $y\in H, y \not\in (\ua G)^{\circ}$. Then, for any open set $U\subseteq X$ with $y\in U$, there exists a $y_U\in U$ such that $y_U\not\in \ua G$. Set $J=\{U\in \mathcal{O}(X): y\in U\}$ ordered as: $\forall U,V\in J$,$U\leq V\Leftrightarrow V\subseteq U.$
Then $J$ is directed.  Easily, one can see that  $(y_U)_{U\in J}$ is a net converging to $y$. Since $G \ll_n y$, there exists $U_0\in J$ such that $x \sqsubseteq y_{U_0}$ for some $x \in G$ by definition. It is a contradiction for $y_{U_0}\not\in \ua G$.  Hence, $y \in (\ua G)^{\circ}$. \ $\Box$

\vskip 3mm
From Lemma \ref{d-quasicontinuity} and Lemma \ref{n-quasicontinuity}, we can see that the notion of $d$-approximation and $n$-approximation coincide for when $X$ is a $d$-quasicontinuous directed space. Since the $n$-approximation does not depend on directed space, a notion of $n$-quasicontinuity for a $T_0$ space can be defined as follows.

\begin{defn}\rm 
A $T_{0}$ topological space $X$ is called $n$-quasicontinuous if for each $x\in X$, there is a directed family $D(x) \subseteq fin_{n}(x) = \{F: F\ \text{is finite}, F \ll_n x\}$ such that $D(x) \rightarrow x$.
\end{defn}

\begin{defn}\rm
A $T_{0}$ topological space $X$ is called locally hypercompact, if for any open subsets $U$ and $x \in U$, there exists some $F \subseteq_{f}  X$ such that $x \in (\ua F)^{\circ} \subseteq \ua F \subseteq U$. $X$ is called hypercompactly based if for any open subsets $U$ and $x \in U$, there exists some $F \subseteq_{f}  X$ such that $x \in (\ua F)^{\circ} = \ua F \subseteq U$.
\end{defn}

Since $G \ll_n H$ iff $H \subseteq (\ua G)^\circ$, definitions of $n$-quasicontinuous and locally compact are equivalent. When $X$ is quasicontinuous, $G \ll_n H$ is also equivalent to $G \ll_d H$. We drop the subscripts $d,n$ when  no confusion can arise.
Similarly, we can define $d$-quasialgebraic space and $n$-quasialgebraic spaces. And the following statements hold. We omit the proof for quasialgebraic cases since it is similar to the proof for quasicontinuous spaces \cite{FK2017, LAW85} and quasialgebraic domains.

\begin{thm}  \rm \cite{FK2017,LAW85} \label{n-quasicontinuous space}
For a $T_{0}$ topological space $X$, the following are equivalent:
\begin{enumerate}
\item[(1)] $X$ is a $d$-quasicontinuous directed space;

\item[(2)] $X$ is  $n$-quasicontinuous;

\item[(3)] $X$ is a locally hypercompact space;
\item[(4)] $\mathcal{O}(X)$ is hypercontinuous.

\end{enumerate}
\end{thm}

\vskip 2mm

\begin{thm}\rm \label{quasialgebraic} 
For a $T_{0}$ topological space $X$, the following are equivalent:
\begin{enumerate}
\item[(1)] $X$ is a $d$-quasialgebraic directed space; 

\item[(2)] $X$ is  $n$-quasialgebraic;

\item[(3)] $X$ is a hypercompactly based space;
\item[(4)] $\mathcal{O}(X)$ is hyperalgebraic;

\end{enumerate}

\end{thm}

\begin{defn}\rm 
A $T_{0}$ topological space is called quasicontinuous if one of the four equivalent conditions in Theorem \ref{n-quasicontinuous space} holds; called quasialgebraic if one of four equivalent conditions in \ref{quasialgebraic} holds.
\end{defn}

Next, we show that for a directed space $X$, $\mathcal{O}(X)$ is algebraic iff it is hyperalgebraic.

\begin{lem}\rm \cite{CK2022} \label{hyperbelow}
Let $X$ be a topological space and $U,V \in \mathcal{O}(X)$. $U \prec V$ iff $U \subseteq \ua F \subseteq V$ for some $F \subseteq_{f} X$.
\end{lem}

\begin{prop} \rm \label{com=hypercom}
Let $X$ be a directed space. $\mathcal{O}(X)$ is algebraic iff it is hyperalgebraic.
\end{prop}
\noindent{\bf Proof.} We need only to show that any compact element of $\mathcal{O}(X)$ is hypercompact. Suppose that $K$ is a compact element of $\mathcal{O}(X)$. Given any $x \in K$, since $K$ is compact, there must be a minimal element $y\in K$ such that $ y \leq x$. Set $A$ to be the set of all minimal elements of $K$, then $K = \ua A$, $A$ is also a compact subset of $X$. Let $B$ be any subset of $A$, we show that $\ua\! A \backslash B$ is an open subset of $X$. Let $D$ be a directed subset of $X$, $x \in\  \ua\! A \backslash B$ and $D \to x$.  By that $\ua A$ is open, there exists some $d\in D$ such that $d \in A$. Since $D$ is directed there are two cases. 

 Case(1): $d$ is the largest element of $D$. Then by $D \to x$, we have $x \leq d$, that is $d \in\ \ua\! A \backslash B$.

 Case(2): there exists some $e \in D$ such that $d < e$. Then $e$ must be in $ \ua\! A \backslash B$.

Thus, $D \cap (\ua\!A\backslash B) \not = \emptyset$, $\ua\! A \backslash B$ is open. Suppose that $A$ has infinite elements. Let $\omega = \{\omega_i\}_{i\in \mathbb{N}} \subseteq A$ and $B_{i} = \ua\! A \backslash \{\omega_j\}_{i \leq j}$ for all $i \in \mathbb{N}$. Then $A = \bigcup_{i \in \mathbb{N}} B_{i}$. By the compactness of $A$, there must exists some $i$ such that $B_{i} =A$, a contradiction. Therefore, $A$ contains only finite elements and then $K= \ua A$ is hypercompact by Lemma \ref{hyperbelow}.

\begin{rem} \rm
Lemma \ref{com=hypercom} does not always hold for general $T_{0}$ spaces. A simple example is $\mathbb{N}$, the set of natural numbers,  endowed with the cofinite topology. Any open subset of it is compact, but not hypercompact.
\end{rem}

By the above discussion, we conclude the relationships between the continuity of a $T_0$ space and the lattice of its open subsets as follows.

\begin{thm}\rm \label{open lattice}
Let $X$ be a $T_0$ space. 
\begin{enumerate}
\item[(1)] $X$ is continuous (algebraic) iff $\mathcal{O}(X)$ is a completely distributive lattice (completely distributive algebraic lattice).
\item[(2)] $X$ is quasicontinuous iff the lattice $\mathcal{O}(X)$ is a hypercontinuous lattice.
\item[(3)] $X$ is quasialgebraic iff $\mathcal{O}(X)$ is a hyperalgebraic lattice. Moreover, if $X$ is a directed space, then $X$ is quasialgebraic iff $\mathcal{O}(X)$ is an algebraic lattice.

\end{enumerate}
\end{thm}

Given any $T_0$ space $X$, denote $C(X)$ the lattice of its closed subsets. We show that the continuity of $X$ can be characterized by the continuity of $C(X)$, and vice versa. We first recall some equivalent conditions for continuity of lattices.

\begin{lem} \rm \cite{CON03} \label{completely distributive}
Let $L$ be a complete lattice. The following statements are equivalent:
\begin{enumerate}
\item[(1)] $L$ is completely distributive;
\item[(2)] $L^{\rm op}$ is completely distributive;
\item[(3)] $L$ is continuous and every element is the sup of co-primes.
\end{enumerate}
\end{lem}

\begin{lem} \rm \cite{CON03,LAW81,VEN90} \label{hyper op} \label{dis hyper}The following statements hold.
\begin{enumerate}

\item[(1)] Let $L$ be a distributive continuous (algebraic) lattice, then $L$ is hypercontinuous (hyperalgebraic) iff $L^{\rm op}$ is hypercontinuous (quasialgebraic).

\item[(2)] A complete lattice $L$ is a quasicontinuous (quasialgebraic) lattice iff $\omega(L)$ is a continuous (algebraic) lattice.
\end{enumerate}

\end{lem}

Now, it is time to show the following theorems.

\begin{prop} \rm \label{quasi op}
Let $X$ be a $T_0$ space. If $C(X)$ is a quasicontinuous lattice, then $X$ is quasicontinuous.
\end{prop}
\noindent{\bf Proof.}
Let $C(X)$ be a quasicontinuous domain. Then for any $x \in X$, there exists some $\{A_1,\dots,A_n\} \ll\ \da x$ in $C(X)$. 
Given any $a_i \in A_i\ (1 \leq i \leq n)$, we have $\{\da a_1,\dots,\da a_n\} \ll\ \da x$.

{\bf Claim 1.} $\{a_1,\dots,a_n\} \ll_n x$ in $X$.

Given any net $(x_i)_I \to x$, $x \in  \overline{\bigcup_{i \in I} \{x_i\} }$. Then $\da x \subseteq \ \overline{\bigcup_{i \in I} \{x_i\} }$. Let $\mathcal{F} = \{\bigcup_{i \in F}\! \da x_i : F \subseteq_f I \}$. Then $\mathcal{F}$ is directed and $ \da x \leq \bigvee^\ua \mathcal{F} = \overline{\bigcup_{i \in I} \{x_i\} } $ in $C(X)$.  There exist some $F\subseteq_f I$ and $a_j$ such that $\da a_j \subseteq \bigcup_{i\in F} \da x_i$. Thus, there exists a $i \in F$ such that $a_j \leq x_i$, i.e., $\{a_1,\dots,a_n\} \ll_n x$ in $X$.

{\bf Claim 2.} For any $x$ and any open neighbourhood $U$ of $x$, there exists some $F \subseteq_f U$ such that $x \in (\ua F)^\circ$.

 Given any directed family $\mathcal{A}=\{A_i\}_I$ of closed subsets of $X$ such that $\da\! x \leq \bigvee \mathcal{A} = \overline{\bigcup_{i \in I} A_i}$, then $U \cap (\bigcup_{i \in I} A_i) \not = \emptyset$. There exists some $i \in I$ such that $A_i \cap U \not = \emptyset$. Pick $b \in A_i \cap U$, then $\da\! b \subseteq A_i$. We have $\{\da\! b: b \in U\} \ll\ \da\! x $ in $C(X)$. Since $C(X)$ is quasicontinuous, there exists a $\{B_1,\dots,B_m\} \subseteq C(X)$ such that $\{\da b : b \in U \} \ll \{B_1,\dots,B_m\} \ll\ \da x $ in $C(X)$. Then  
 $$\forall 1 \leq i \leq m,\ \{\da b: b \in U\} \ll B_i = \overline{\bigcup_{F\subseteq_f B_i} \da F},$$ 
 which implies that there exists some $F \subseteq_f B_i$ and $b_i \in U$ such that $\da b_i \subseteq \da F \subseteq B_i$. Thus, $b_i \in B_i \cap U$. It follows that $\{\da b_1,\dots,\da b_m\} \ll\ \da x$ in $C(X)$. Then $\{b_1,\dots,b_n\} \ll_n x$ by Claim 1, i.e., $x \in (\ua\{b_1,\dots,b_n\})^\circ$. $\Box$

\begin{thm}\rm
Let $X$ be a $T_0$ space.
\begin{enumerate}
\item[(1)] $X$ is continuous iff the lattice $C(X)$ is a continuous lattice iff $C(X)$ is completely distributive.
\item[(2)] $X$ is quasicontinuous iff the lattice $C(X)$ is a quasicontinuous lattice.
\item[(3)] $X$ is algebraic iff the lattice $C(X)$ is an algebraic lattice.
\item[(4)] $X$ is quasialgebraic iff the lattice $C(X)$ is a quasialgebraic lattice.
\end{enumerate}
\end{thm}
\noindent{\bf Proof.}
(1) Let $C(X)$ be a continuous domain. Given any $x \in X$, $\da x$ is co-prime in $C(X)$. For any $F \in C(X)$, $F = \bigcup \{\da x : x \in F\} = \bigvee  \{\da x : x \in F\}$. Thus, $C(X)$ has enough co-primes. It follows that $C(X)$ and $\mathcal{O}(X)$ are completely distributive lattices by Theorem \ref{completely distributive}. Then,  $X$ is a continuous space by \ref{open lattice}. Conversely, let $X$ be a continuous space. Then $\mathcal{O}(X)$ is completely distributive. It follows that $C(X)$ is completely distributive and hence continuous.

(2) By Proposition \ref{quasi op}, if $C(X)$ is quasicontinuous, then $X$ is quasicontinuous. Conversely, let $X$ be quasicontinuous. Then $\mathcal{O} (X)$ is hypercontinuous by Theorem \ref{open lattice}. From Theorem 
\ref{hyper op}, we know that $C(X)$ is quasicontinuous.

(3) It is analogous to the proof for (1), since the completely distributive algebraic lattice is selfdual \cite{Erne91}.

(4) Let $C(X)$ be a quasialgebraic lattice. Then $\omega (C(X)) = \upsilon( \mathcal{O}(X))$ is algebraic by Lemma \ref{hyper op} and $\mathcal{O}(X)$ is hypercontinuous by (2) and Theorem \ref{open lattice}. Then $\sigma(\mathcal{O}(X)) = \upsilon(\mathcal{O}(X))$ is algebraic and completely distributive. Thus, $\mathcal{O}(X)$ is algebraic and then hyperalgebraic. Conversely, let $X$ be a quasialgebraic space. Then $\mathcal{O}(X)$ is distributive hyperalgebraic lattice. $C(X)$ is a quasialgebraic lattice by Theorem \ref{dis hyper}.

\section{Core compactness of directed spaces}

Core compactness can be viewed as a weaker continuity property than quasicontinuity. Not all core compact spaces are directed spaces. All nontrivial compact $T_{2}$ spaces are the examples and Example \ref{exand}(1) is a locally compact sober space, which is neither a directed space nor a $T_{1}$ space. In \cite{LAW85}, Lawson gave some equivalent conditions for a $T_{0}$ topological space $X$ to be quasicontinuous. One of the key equivalent conditions is that for any $T_{0}$ topological space $Y$, $X \times Y =X \bigotimes Y$, where $X \bigotimes Y$ is the tensor product of $X$ and $Y$. In \cite{CON03}, some equivalent conditions for a $T_{0}$ topological space to be core compact was given.  Denote $\Sigma P$ the topological space $(P, \sigma(P))$. 
In particular, it was proved that for a dcpo $P$, $P$ is core compact iff for any dcpo $Q$, $\Sigma(P \times Q) = \Sigma P \times \Sigma Q$.

In this section, we investigate the core compactness of a directed space. We show that for any two directed spaces $X$ and $Y$, their tensor product $X \bigotimes Y$ is same to their categorical product $X \otimes Y$ in $\bf DTop$. Similar to Scott spaces, a directed space is core compact iff for any directed spaces $Y$, $X \times Y = X \otimes Y$. 
Finally, we give more equivalent conditions for a directed space to be core compact.

\begin{defn} \rm \cite{CJ1970}
For any two topological spaces $X,Y$, the tensor product $X \bigotimes Y$ of $X$ and $Y$ has the same carrier set of $X \times Y$. A set $W$ is open in $X \bigotimes Y$ if for all $(x,y) \in X \times Y$, the slices $W_{x} =\{y \in Y: (x,y) \in W\}$ and $W^{y} = \{x \in X : (x,y) \in W\}$ are open in $Y$ and $X$ respectively.
\end{defn}

Given any two topological spaces $X,Y$, a map from a topological space $(X \times Y, \tau)  $ to $Z$ is called separately continuous if it is continuous at each argument, i.e., for any $(x_{0},y_{0}) \in X\times Y$, the maps $f_{x_{0}}: Y \to Z$ and $f^{y_{0}}: X \to Z$ are continuous, where $f_{x_{0}}(y) = f(x_{0},y), f^{y_{0}}(x) = f(x,y_{0})$. The topology of tensor product is also called the topology of separate continuity \cite{LAW85}, which is the weak topology determined by all separately continuous functions from the product. A map from $X \bigotimes Y$ is continuous iff it is separately continuous \cite{CJ1970}.

\begin{lem}  \label{tensor=product} \rm
For any two directed spaces $X$ and $Y$, $X \bigotimes Y = X \otimes Y$.
\end{lem}

\noindent{\bf Proof.}
Assume that $W$ is open in $X \otimes Y$. Given any $(x,y) \in W$, we show that the slice $W^{y}$ is open in $X$. For any directed subset $(x_{i})_{I}$ of $X$ with $(x_{i})_{I} \to x \in W^{y}$ in $X$, since $((x_{i},y))_{I} \to (x,y)$ in $X \otimes Y$, then $((x_{i},y))_{I}$ is finally in $W$. It follows that $(x_{i})_{I}$ is finally in $W^{y}$. So, $W^{y}$ is open in $X$. Similarly, $W_{x}$ is open in $Y$. Therefore, $W$ is open in $X \bigotimes Y$.

Conversely, assume that $W$ is open in $X \bigotimes Y$.
Given any directed subset $(x_{i})_{I} \to x$ in $X$, then for $D = ((x_{i},y))_{I}$ we have $D \to (x,y) \in W$ in $X \times Y$. Since $W^{y}$ is open in $X$ and $X$ is a directed space,  there exists some $i_{0} \in I$ such that $x_{i_{0}}\in W^{y}$. Then $(x_{i_{0}},y) \in W \cap D$. Similarly, for any directed subsets $((x,y_{i}))_{I} \to (x,y) \in W$, there exists some $y_{i}$ such that $(x,y_{i}) \in W \cap D$. Therefore, $W$ is open in $X \otimes Y$ by Lemma \ref{directed product}. $\Box$

\begin{cor}\rm Let $X,Y,Z$ be directed spaces. A map $f: X \otimes Y \to Z$ is continuous iff it is separately continuous.
\end{cor}

We recall the following condition for a $T_{0}$ topological space to be core compact.

\begin{thm}\rm \cite{CON03} \label{core compact of scott space} 
Let $X$ be a $T_{0}$ space. Then the following statements are equivalent.
\begin{enumerate}[(1)]
 \item $ \mathcal{O}(X) $ is a continuous lattice.
 \item The set $\{ (U,x) \in \mathcal{O}(X) \times X : x \in U \}$ is open in $  \Sigma \mathcal{O}(X) \times X $.

\end{enumerate}
\end{thm}

\begin{thm} \label{core compact vs product} \rm
Let $X$ be a directed space. Then the following conditions are equivalent.
\begin{enumerate}[(1)]
\item $X$ is core compact.  
\item $X \otimes Y = X \times Y$ for any directed space $Y$.
\item $\Sigma \mathcal{O}(X) \times X  =  \Sigma \mathcal{O}(X) \otimes X$.
\end{enumerate}
\end{thm}

\noindent{\bf Proof.}
$(3) \Rightarrow (1)$. Assume that $ \Sigma \mathcal{O}(X) \times X  =  \Sigma \mathcal{O}(X) \otimes X$. $\Sigma \mathcal{O}(X)$ is a directed space. To show that $X$ is core compact, we need only to show that $E = \{ (U,x) \in\mathcal{O}(X) \times X : x \in U \}$ is open in $\Sigma \mathcal{O}(X) \times X$ by Lemma \ref{core compact of scott space}. Then it is equivalent to showing that $E$ is open in $\Sigma \mathcal{O}(X) \otimes X$. 
Assume that a directed set $((U_{i},x_{i}))_{i \in I}$ converges to $(U,x) \in E$ in $\Sigma \mathcal{O}(X) \times X$. Then $(x_{i})_{I} \to x$, $(U_{i})_I \to U$. Therefore, there exists some $i_{0} \in I$ such that $x_{i_{0}} \in U$ and $U \subseteq \bigcup_{i \in I} U_{i}$. Then $\exists i_{1} \in I$ such that $x_{i_{0}} \in U_{i_{1}}$. Let $i_{0}, i_{1} \leq i_{2} $, we have $(U_{i_{2}},x_{i_{2}}) \in E$.

$(1) \Rightarrow (2)$. Suppose that $X$ is core compact. We show that $X \otimes Y = X \times Y$ for any directed space $Y$. 
Since $X \otimes Y$ is finer than $X \times Y$, we need only to show that every open subset $U$ of $X \otimes Y$ is open in $X \times Y$. Given $(x_{0},y_{0}) \in U$, let $V = \{x\in X: (x,y_{0}) \in U\}$. By Lemma \ref{tensor=product}, $V$ is an open subset of $X$. Since $X$ is core compact, there exists a family of open subsets $\{V_{n}: n \in N \}$ such that $$ x_{0} \in V_{0} \ll  \cdots \ll V_{n+1} \ll V_{n} \ll \cdots \ll V_{2} \ll V_{1} \ll  V.$$

Let $W = \bigcup_{1 \leq n}\{y\in Y: V_{n} \times \{y\} \subseteq U \}$. Obviously, $y_{0} \in W$.
Since $V_{0} \subseteq V_{n}$ for any $n \geq 1$, then $V_{0} \times W \subseteq U$. We need only to show that $W$ is an open subset of $Y$.
Given any directed subset $D \to y \in W $ in $Y$, there exists some $n$ such that $V_{n} \times \{y\} \subseteq U$. For any $x \in V_{n}$, $((x,d))_{d \in D} \to (x,y) \in U$ in $X \otimes Y$. 
Thus, there is some $d_{x} \in D$ such that $(x,d_{x}) \in U$. Then there exists an open neighborhood $V_{x}$ of $x$ such that $V_{x} \times \{d_{x}\} \subseteq U$. Since $V_{n+1} \ll V_{n} \subseteq \bigcup_{x\in V_{n}} V_{x}$, we have that $V_{n+1} \subseteq \bigcup_{i=1}^{n} V_{x_{i}}$ for some finite subset $B = \{x_{1},\cdots, x_{n} \}$ of $V_{n}$. Since $B$ is finite, there exists some $d_{0} \in D$ such that $d_{x_{i}} 
\leq d_{0}$. Then $V_{n+1} \subseteq \{x\in X: (x,d_{0}) \in U \}$ and thus $V_{n+1} \times \{d_{0}\} \subseteq U$, i.e., $d_{0} \in W$. Therefore, $W$ is open in $Y$. Then $(x_0,y_0) \in V_0 \times W \subseteq U$, i.e., $U$ is open in $X \times Y$.

$(2)\Rightarrow (3)$. Obviously.
$\Box$

\vskip 2mm

Now, we give some more equivalent conditions for a directed space to be core compact. We denote $\Sigma 2$ to be the Sierpinski space, i.e., the set $\{0,1\}$ endowed with the topology $\{\emptyset, \{1\}, \{0,1\}\}$.

\begin{lem} \rm \cite{LAW85}  \label{pointwise topology vs upper}
The topology of pointwise convergence on $[X \to \Sigma 2]_p$ is the upper topology, which corresponds to the upper topology on $\mathcal{O}(X)$.
\end{lem}

\begin{prop}\rm \cite{CON03}\label{Isbell} 
If $X$ is a space such that $\mathcal{O}(X)$ is a continuous lattice and $Y$ is an injective space, then $[X \to Y]_I$ is injective. In particular, the Isbell topology on $TOP(X,Y)$ is the Scott topology.
\end{prop}

\begin{thm} \rm
Let $X$ be a directed space. The following conditions are equivalent.

\begin{enumerate}
	\item $[X \to Y]$ is injective for all injective $T_{0}$ spaces $Y$.
	\item $[X \to \Sigma 2]$ is injective.
	\item $\mathcal{O}(X)$ is continuous.
	\item $\{(U,x): x \in U \}$ is open in $\Sigma \mathcal{O}(X) \times X $.
	\item The evaluation map $ev : [X \to \Sigma 2] \times X \to \Sigma 2$ is continuous.
	\item For all directed spaces $Y$, $X \bigotimes Y = X \otimes Y = X \times Y$.
	\item For all directed spaces $Y,Z$, if a map $f: X \times Y \to Z $ is separately continuous, then it is jointly continuous.
	\item For any directed space $Y$, the evaluation map $[X \to Y] \times X \to Y$ is continuous.

	\item The natural map $[Z \times X \to Y ] \to [Z \to [X \to Y]]$ is onto (and a homeomorphism) for all directed spaces $Y$ and $Z$.

\end{enumerate}

\end{thm}

\noindent{\bf Proof.}
$(1) \Rightarrow (2)$. Obviously.

$(2) \Leftrightarrow (3)$.  Since $[X \to \Sigma 2] = \mathcal{D}([X \to \Sigma 2 ]_p) = \Sigma \mathcal{O}(X)$ and an injective space is a continuous lattice endowed with the Scott topology, then $\mathcal{O}(X)$ is continuous iff $[X \to\Sigma 2]$ is injective.

$(3) \Rightarrow (1)$. Assume that $Y$ is an injective space. By Proposition \ref{Isbell}, $[X \to Y]_I = \Sigma(TOP(X,Y))$ is an injective space. Thus, $TOP(X,Y)$ is a continuous lattice. Since the topology of pointwise convergence is coarser than the Isbell topology, $[X \to Y]_p$ is coarser than $\Sigma(TOP(X,Y))$ and then $[X \to Y] = \mathcal{D} ([X \to Y]_p) = \Sigma(TOP(X,Y)) =  [X \to Y]_I$.

$(3) \Leftrightarrow (4)$. By Theorem \ref{core compact of scott space}.

$(4) \Leftrightarrow (5)$. It is a direct conclusion by the fact that $[X \to \Sigma 2] =  \Sigma \mathcal{O}(X)$. 

$(3) \Leftrightarrow (6)$. By Theorem \ref{core compact vs product} and Lemma \ref{tensor=product}.

$(6) \Rightarrow (7)$. $f$ is separately continuous is equivalent to that $f$ is continuous from $X \bigotimes Y $ to $Z$. Thus $f$ is jointly continuous.

$(7) \Rightarrow (6)$. Given any directed space $Y$, let $Z = X \otimes Y$. Then $Z$ is also a directed space. The identity map $id : X \times Y \to X \otimes Y$ is separately continuous since $X \otimes Y = X \bigotimes Y$. Thus, $id : X \times Y \to X \otimes Y$ is continuous. Then $X \times Y = X \otimes Y$.

$(7) \Rightarrow (8) \Rightarrow (5)$. Since $ev: [X \to Y] \otimes X \to Y$ is continuous, i.e., $ev: [X \to Y] \times X \to Y$ is separately continuous, then $ev$ is continuous from $[X \to Y] \times X$. $(8) \Rightarrow (5)$ is obvious.

$(8) \Rightarrow (9)$. Since the natural map $[Z \otimes X \to Y] \to [Z \to [X \to Y]]$ is a homeomorphism for all directed spaces $Y$ and $Z$ (see \cite{L2017,YK2015}), we need only to show (6). This has been proved.

$(9) \Rightarrow (8)$. Let $Z = [X \to Y]$. Then the inverse of identity map $id : [[X \to Y] \to [X \to Y]]$ is $ev: [[X \to Y] \times X \to Y]$. $\Box$

\section{Questions related to core compactness}

It is well known that the spectrum with hull-kernel topology of a completely distributive lattice (resp. a distributive hypercontinuous lattice) is exactly a continuous (resp. quasicontinuous) dcpo endowed with the Scott topology \cite{LAW81,Hoff1981,Lawson1979}.  In this section, some conclusions for Scott spaces are extended to directed spaces. These conclusions are close related to the long-standing open problem that, for which distributive continuous lattice its spectrum is exactly a sober locally compact Scott space (see \cite[Problem 528]{Mill1990}).

Given any poset $P$, $\upsilon(\Gamma (P)) = \sigma(\Gamma(P))$ is a necessary condition for $\Sigma P$ to be core compact \cite{CK2022}. We show that for any directed space $X$, $X$ is core-compact iff $(C(X), \sigma(C(X)))$ is sober and locally compact with $\sigma(C(X))=\upsilon(C(X)).$

Given a topological space $X$, we denote $\prod\limits^n X$ to be topological product of $n$ copies of $X$. For any $n\in \mathbb{N}$, define a map $s_n : \prod\limits^n X \to  \Sigma(C(X))$ as follows: $\forall (x_1,x_2,\dots,x_n) \in \prod\limits^n X$,
$$s_n(x_1,x_2,\dots,x_n) =\ \da\! {\{x_1,x_2,\dots,x_n\}}.$$

\begin{prop} \rm \label{opop}
For a topological space $X$, $\sigma(C(X)) = \upsilon(C(X))$ holds iff $s_n$ is continuous for all $n \in \mathbb{N}$.

\end{prop}

\noindent{\bf Proof.} 
Assume that $\sigma(C(X)) = \upsilon(C(X))$. Given any $F \in C(X)$, 
$$ s_n^{-1} (\da_{C(X)} F) \ = \{\{x_1,x_2,\dots,x_n\} \in \prod\limits^n X: \overline{\{x_1,x_2,\dots,x_n\}} \subseteq F\} = \prod\limits^n F$$ is a closed subset of $\prod\limits^n X$. Thus,  $s_n: \prod\limits^n X \to (C(X),\upsilon(C(X)))$ is continuous. Then $s_n: \prod\limits^n X \to \Sigma(C(X))$ is continuous.

For the converse, assume that $s_n$ is continuous for all $n \in \mathbb{N}$. Let $\mathcal{U}$ be an open subset of $\Sigma(C(X))$ and $A\in\mathcal{U}$. Without loss of generality, we assume $A\not=\emptyset$. Note that since $A=\bigcup\{\da\! F: F\subseteq_f\!A\}$ and $\{\da F: F\subseteq_fA\}$ is  a directed family in $C(X)$, there exists a non-empty finite subset $F$ of $ A$ such that $\da \! F\in \mathcal{U}$. Let $F=\{x_1,x_2,\ldots,x_n\}$, then $s_n(x_1,x_2,\ldots,x_n)=\da\! F\in \mathcal{U}$. It follows that  $(x_1,x_2,\ldots,x_n)\in s_n^{-1}(\mathcal{U})$. By the continuity of $s_n$, there exists a family of open subsets $U_k (1 \leq k \leq n)$
such that $U_1\times U_2\times\cdots\times U_n$ is open in $\prod\limits^n X$ and 
$$(x_1,x_2,\ldots,x_n)\in U_1\times U_2\times\cdots\times U_n\subseteq s_n^{-1}(\mathcal{U}).$$
Since $x_k\in A$ for $1\leq k\leq n$, we have $A\in \Diamond U_k=\{B\in C(X): \ B\cap U_k\not=\emptyset\}$. It follows that $A\in \bigcap\limits_{k=1}^n\Diamond U_k\in \upsilon(C(X))$.
For any $B\in \bigcap\limits_{k=1}^n\Diamond U_k$, there exists $y_k\in B\cap U_k$ for $1\leq k\leq n$. Since $(y_1,y_2,\dots,y_n)\in s_n^{-1}(\mathcal{U})$, we have $\bigcup\limits_{k=1}^n\da y_k\in \mathcal{U}$. It follows that $B\in \mathcal{U}$, i.e., $A\in\bigcap\limits_{k=1}^n\Diamond U_k\subseteq \mathcal{U}$. $\Box$

\begin{prop}\rm \label{Lop is sober locally compact} \cite{CK2022}
Let $L$ be a continuous lattice, if $L$ satisfies the condition that $\upsilon(L^{op})=\sigma(L^{op})$, then $(L^{op}, \sigma(L^{op}))$ is a sober and locally compact space.
\end{prop}

\begin{prop}\rm \label{core op} Let $X$ be a directed space.\ If $X$ is core compact,\ then $\sigma(C(X)) = v(C (X))$. Moreover, $(C(X), \sigma(C(X)))$ is sober and locally compact.
\end{prop}
\noindent{\bf Proof.}
 Let $X$ be a core compact directed space. We have that for every $n \in \mathbb{N}$, $\mathcal{D}(\prod\limits^n X) = \prod\limits^n X$ by Theorem \ref{core compact vs product}. 
Then $s_n$ is continuous from $ \prod\limits^n X  $ to $ \Sigma(C(X))$ iff it is continuous from $\mathcal{D} (\prod\limits^n X)  $ to $ \Sigma(C(X))$.

We show that $s_n: \mathcal{D}(\prod\limits^n X) \to \Sigma(C(X))$ is continuous, i.e., $s_{n}$ preserves $D \to x$ for every $(D,x) \in DLim(\prod\limits^n X)$. Let $\{(x_{1i},x_{2i},\ldots,x_{ni}): \ i\in I\}$ be a directed subset of  $\prod\limits^n X$ converging to $(x_{1},x_{2},\ldots,x_{n})$ in $\prod\limits^n X$. Then for each $1 \leq k \leq n$, $\{ (x_{ki}): i \in I\}$ converges to $x_{k}$ by the definition of topological product. We have 

$$s_n((x_{1},x_{2},\ldots,x_{n})) = \bigcup \limits_{k=1}^n \downarrow x_{k} \subseteq \overline{\bigcup\limits_{k=1}^n \bigcup_{i\in I} \downarrow x_{ki}}= \overline{\bigcup_{i\in I}\bigcup\limits_{k=1}^n\downarrow x_{ki}} = \bigvee_{i\in I}s_n(x_{1i},x_{2i},\ldots,x_{ni}).$$Thus $s_n$ is a continuous map from $\mathcal{D}(\prod\limits^n X)$ into $\Sigma(C(X))$.

 By Proposition \ref{opop}, we have $\sigma(C(X)) = v(C(X))$. Letting $L = \mathcal{O}(X)$, then $L$ is a continuous lattice and $C(X) = L^{op}$. $(C(X),\sigma(C(X)))$ is sober and locally compact by Proposition \ref{Lop is sober locally compact}. $\Box$

\vskip 2mm

In \cite{CK2022}, an adjunction between  $\sigma(P)$ and $\sigma(\Gamma(P))$ serves as a useful tool in studying the relation between  $P$ and $\Gamma(P)$. It can be extended to directed spaces as well.

\begin{defn} \rm \label{eta}
Given a directed space $X$, we define a map $\eta: X\rightarrow C(X)$ and a map $\Diamond:\mathcal O(X)\rightarrow \sigma(C(X))$ as follows: $\forall x\in X$, $\forall U\in \mathcal{O}(X)$,
$$\eta(x)=\da x, \ \Diamond(U)=\{A\in C(X): A\cap U\not=\emptyset\}.$$
Define $\eta^{-1}:\sigma(C(X))\rightarrow \mathcal{O}(X)$ by $\eta^{-1}(\mathcal U)=\{x\in X:\ \da x\in \mathcal{U}\}$. 
\end{defn}

Then we have the following result.

\begin{prop}\rm \label{ajuction}
For a directed space $X$, both $\eta^{-1}$ and $\Diamond$ preserve arbitrary $\sup$s. Moreover, $\Diamond \circ \eta^{-1}\leq 1_{\sigma(C(X))}$ and $\eta^{-1}\circ \Diamond=1_{\mathcal{O}(X)}$.
\end{prop}

\noindent{\bf Proof.}
 Let $\{U_{i}: i \in I\}$ be any subset of $\sigma(C(X))$. Then $\Diamond(\bigcup_{i\in I}U_i)=\{A\in C(X): A\cap\bigcup_{i\in I}U_i\neq\emptyset\}=\{A\in C(X): \exists i\in I, A\cap U_i\neq\emptyset\} = \bigvee_{i\in I}\Diamond(U_i)$. $\eta$ is the special case of $s_{n}$ for $n =1$. Thus it is continuous. Then $\eta^{-1}$ preserves arbitrary sups. Given any $U\in\ \mathcal{O}(X),\
x\in\eta^{-1}(\Diamond(U))\Leftrightarrow\ \eta(x)\in \Diamond(U)\Leftrightarrow\ \da x\cap U\neq\emptyset\Leftrightarrow x\in U$, hence $\eta^{-1}\circ \Diamond=1_{\mathcal{O}(X)}$. For any $\mathcal{U}\in\sigma(C(X))$, $A\in \Diamond\circ\eta^{-1}(\mathcal U)\Leftrightarrow\ A\cap\eta^{-1}(\mathcal U)\neq\emptyset\Rightarrow A\in\mathcal U$, i.e., $\Diamond\circ\eta^{-1}\leq 1_{\sigma(C(X))}$. $\Box$

\begin{thm}\rm\label{corecom closedsubsets}
Let $X$ be a directed space, then 
$X$ is core-compact
iff  $(C(X), \sigma(C(X)))$ is core compact 
iff $(C(X), \sigma(C(X)))$ is sober and locally compact with $\sigma(C(X))=\upsilon(C(X))$.

\end{thm}

\noindent{\bf Proof.}
Suppose that $X$ is core compact. By Proposition \ref{core op}, $(C(X), \sigma(C(X)))$ is sober and locally compact with $\sigma(C(X))=\upsilon(C(X))$. Conversely, suppose that $(C(X), \sigma(C(X)))$ is core compact, i.e., $\sigma(C(X))$ is continuous. By Proposition \ref{ajuction}, $\mathcal{O}(X)$ is continuous, i.e., $X$ is core compact.
$\Box$

\vskip 2mm

\begin{prob} \rm
Let $X$ be a $T_0$ space and $\Sigma C(X)$ be core compact. Must $X$ be core compact?
\end{prob}

The adjunction $(\eta^{-1},\Diamond)$ seems only holds for directed spaces. A natural question arises that whether a topological space $X$ that makes the map $\eta$ in Definition \ref{eta} continuous is a directed space? When $L$ is a complete lattice endowed with a topology coarser than $\sigma(L)$, the answer is positive. However, for other cases, we still do not know the answer.

\begin{lem} \rm \label{retract}
Any retract of a directed space is a directed space.
\end{lem}
\noindent{\bf Proof.} Let $X$ be a topological space and $Y$ be a directed space. Suppose that $i:X \to Y$ and $r:Y \to X$ are continuous maps and $r \circ i = id_X$. We need only to check that any directed open subset $U$ of $X$ is open in $X$. Noticing that $U=(r\circ i)^{-1}(U) = i^{-1}(r^{-1}(U))$, we need only to show that $r^{-1}(U)$ is open in $Y$. Given any directed subset $D$ of $Y$ and $D \to y \in r^{-1}(U)$, then $r(D) \to r(y) \in U$. There exists some $d \in D$ such that $r(d) \in U$, i.e., there exists some $d \in r^{-1} (U)$. Thus, $r^{-1}(U)$ is open in $Y$. $\Box$

\begin{prop} \rm
Let $L$ be a complete lattice and $X=(L,\tau)$ with $\tau \subseteq \sigma(L)$. If $\eta : X \to \Sigma(C(X))$ is continuous, then $X$ is a Scott space.
\end{prop}
\noindent{\bf Proof.} Since $L$ is a complete lattice, we have $\da (\inf A) = \bigcap_{x_{i} \in A} \da x_{i}$, that is, $\eta:L \to C(X)$ preserves all infs. Then, there exists a right adjoint $d:C(X) \to L$ such that $d(F) = \inf \eta^{-1}(\ua_{C(X)} F)= \sup F$. Thus, $(\eta,\sup)$ form a pair of adjunction and $\sup:C(X) \to L$ preserves all sups. Then the map $\sup:\Sigma(C(X)) \to X$ is continuous. It is easy to check that $\sup \circ \eta = id_X$. Thus $X$ is a retract of $\Sigma(C(X))$ and a directed space by Lemma \ref{retract}. Then $X$ is a Scott space since the Scott topology is the coarsest topology on $L$ such that it is a directed space. $\Box$

\vskip 3mm

By Theorem \ref{corecom closedsubsets}, if the spectrum space of a distributive continuous lattice $L$ is a directed space, then $\sigma(L^{op})=\upsilon(L^{op})$ must hold. 
Particularly, the reverse holds when $L$ is algebraic \cite{CK2022,ERNE2009}. So, we emphasize the following open question:

\begin{prob}\rm \label{P1}
Is the hull-kernel topology of the spectrum ${\mathrm{Spec} L}$ equal to the Scott topology when $L$ is a distributive continuous lattice with $\sigma(L^{op})=\upsilon(L^{op})$? 

Equivalently, let $X$ be a sober and core compact space with $\upsilon(C(X))=\sigma(C(X))$. Is $X$ a directed space?

\end{prob}

An alternately relative problem is :
\begin{prob}\rm \label{P2}
Is the soberification of a core compact directed space a directed space (Scott space)?

\end{prob}

Obviously, so Problem \ref{P2} must be if Problem \ref{P1} is  affirmative. There exists a non-continuous spatial complete lattice $L$ with $\sigma(L^{op})=\upsilon(L^{op})$, but its spectrum is not a Scott space.

\begin{exmp} \rm
Let $\mathbb{J}$ be the classical non-sober dcpo given by Johnstone \cite{Stone}. Denote $\mathbb{J} = \mathbb{N} \times (\mathbb{N} \cup \{\omega\})$.  $(m_1,n_1) \leq (m_2,n_2)$ iff $m_1 =m_2, n_1 \leq n_2$ in $\mathbb{N}$ or $n_2 = \omega, n_1 \leq m_2$ or $m_1 =m_2, m_2 = \omega$.

It satisfies that $\sigma(\Gamma(\mathbb{J}))=\upsilon(\Gamma(\mathbb{J})) $. Set $L=\sigma(\mathbb{J})$. The spectrum of $L$ is just to add a top element to $\mathbb{J}$, i.e., ${\rm Spec}L=\mathbb{J}\cup\{\top\}$, which is not sober with its Scott topology. Hence, the hull-kernel topology of ${\rm Spec}L$ is not equal to the Scott topology. It is also an example of soberification of a directed space is not a directed space.

(1) 
$\sigma(\Gamma(\mathbb{J})) = v(\Gamma(\mathbb{J}))$. Given any closed subset $\mathcal{A}$ of $\Sigma(\Gamma(\mathbb{J}))$, $\bigcup \mathcal{A}$ is a lower subset of $\mathbb{J}$.
We show that $\bigcup \mathcal{A}$ is closed in $\Sigma\mathbb{J}$. Given any directed subset $D$ in $\bigcup \mathcal{A}$, either $D$ contains a largest element $x$ of $D$ or cofinal with one chain $\{m\} \times \mathbb{N}$ of $\mathbb{J}$ and has a maximal element $(m,\omega)$ of $\mathbb{J}$ as the supremum. For the first case, $\sup D = x \in  \bigcup \mathcal{A}$; for the second case, since $\mathcal{A}$ is a Scott closed subset of $\Gamma(\mathbb{J})$, $\{\da d: d \in D\} \subseteq \mathcal{A}$. Then $\da (m,\omega)$ must be in $\mathcal{A}$ and $(m,\omega) \in \bigcup \mathcal{A}$. Thus, $\bigcup \mathcal{A}$ is Scott closed. Then let $A = \eta^{-1}(\mathcal{A}) = \{x \in \mathbb{J} :\ \da x \in \mathcal{A}\} = \bigcup \mathcal{A}$. $A$ is closed in $\Sigma(\mathbb{J})$. If $A$ contains infinite maximal points of $\mathbb{J}$, then for each $(m,\omega) \in A$, let $B_m = \mathbb{N} \times \{1,2,\dots,m\}$. Then $\mathcal{B} = \{B_m: (m ,\omega) \in A\} \subseteq \mathcal{A}$ form a directed subset of $\Gamma(\mathbb{J})$ and $\bigvee \mathcal{B} = \overline{\bigcup \mathcal{B}} = \mathbb{J}$. Thus, $\mathbb{J} \in \mathcal{A}$, $\mathcal{A} = \Gamma(\mathbb{J})$. 

Now we consider that $A$ contains only finite maximal points of $\mathbb{J}$. It is easy to see that the topology on $A$ inherent from $\Sigma \mathbb{J}$ is equal to the Scott topology. Then, given any open subset $U$ of $\Sigma A$ and $x \in U$, there must be a compact open subset $K$ such that $x \in K \subseteq U$. We need only to let $K =\ \ua\! x\ \cup \ua\! \{x_1, \dots, x_m\}$, where $x_i (1 \leq i \leq m)$ is a picked element that is lower than each maximal element $(n_i, \omega) \in\ \ua\! x \cap U \cap A$. Then $K$ is compact and open. So, $(A,\sigma(A))$ is a locally compact space and hence a core compact space. By Proposition \ref{core op}, $\upsilon(\Gamma(A)) = \sigma(\Gamma(A))$. Since $\mathcal{A}$ is a closed subset of $\Sigma(\Gamma(\mathbb{J}))$, $\mathcal{A}$ is closed in $\Sigma(\Gamma(A))$ and then closed in
$(\Gamma(A),\upsilon(\Gamma(A))$. Thus, $\mathcal{A}$ can be denote as an intersection of a family of finitely generated lower sets in $\Gamma(A)$, and then also an intersection of a family of finitely generated lower sets in $\Gamma(\mathbb{J})$. Thus, $\sigma(\Gamma(\mathbb{J})) = v(\Gamma(\mathbb{J}))$.

(2) ${\rm Spec}L = \mathbb{J} \cup \{\top\}$. It is easy to see that a closed subset of $\Sigma \mathbb{J}$ either contains finite maximal elements, or is equal to the whole space. For the first case, it is not irreducible. For the second case, it is irreducible since any closed subset that contains infinite maximal points of $\mathbb{J}$ must be equal to $\mathbb{J}$. Thus, the only non-trivial irreducible  closed subset of $\Sigma \mathbb{J}$ is $\mathbb{J}$. Then ${\rm Spec}L$ is order isomorphic to  $\mathbb{J} \cup \{\top\}$.

(3) Hull-kernel topology of ${\rm Spec}L$ is not equal to the Scott topology. By definition, a nonempty set $U$ is an open set of the hull-kernel topology iff $U = V \cup \{\top\}$, where $V$ is a non-empty open set of $\Sigma \mathbb{J}$. There is no directed subset of $\mathbb{J}$ in ${\rm Spec}L $ whose supremum is $\top$. Thus, $\{\top\}$ is Scott open in ${\rm Spec}L$. So, the two topologies are not equal.

\end{exmp}

\begin{thm}\rm \cite{CK2022,ERNE2009}\label{hyper alge}
Let $L$ be a continuous lattice. Consider the following conditions:
\begin{enumerate}
\item[(1)] $\sigma(L)=\upsilon(L)$,
\item[(2)] $\sigma(L^{op})=\upsilon(L^{op})$,
\item[(3)] every upper set closed in the dual Scott topology $\sigma(L^{op})$ is compact in the Scott $\sigma(L)$,
\item[(4)] the hull-kernel topology of the spectrum ${\rm Spec}L $ is equal to its Scott topology.
\end{enumerate}
Then $(1)\Rightarrow (2)\Leftrightarrow(3)$. When $L$ is distributive, one has $(4)\Rightarrow (2)$. Additional, if $L$ is distributive and algebraic, then $(1)\Leftrightarrow (2) \Leftrightarrow (3)\Leftrightarrow(4)$.
\end{thm}

In Theorem \ref{hyper alge}, condition(1) is equivalent to $L$ being hypercontinuous. J. Lawson  gave an important example of a meet-continuous non-continuous lattice $W$ such that  the  Scott topology $\sigma(W)$ is continuous (see \cite[Theorem VI-4.5]{CON03}). Let $L=\sigma(W)$. Then $\sigma(L^{op})=\upsilon(L^{op})$. However, since $W$ is not quasicontinuous, it follows that $\sigma(L)\not=\upsilon(L)$.

\begin{lem} \rm \cite{WEN18}\label{closed compact op}
Let $L$ be a complete lattice, $F \subseteq L$. $F$ is closed in $(L,\sigma(L))$ iff it is compact saturated in $(L^{op},\upsilon(L^{op}))$.
\end{lem}

\begin{lem} \rm \cite{CK2022,XY2017}\label{up sober}\label{scott filtered}
Let $L$ be a complete lattice. Then $(L,\upsilon(L))$ is sober and $(L,\sigma(L))$ is well-filtered. 
\end{lem}

\begin{lem} \rm \cite{XY2020}\label{core sober}
A topological space $X$ is core compact and well-filtered iff $X$ is locally compact and sober.
\end{lem} \rm

A complete lattice $L$ is said to be lean if the condition (3) of Theorem \ref{hyper alge} holds \cite{Huth2000}. We finally give an equivalent condition for $\sigma(L^{op}) = \upsilon(L^{op})$.

\begin{thm}\rm \label{lean}
For a continuous lattice $L$, the following two conditions are equivalent to each other:
\begin{enumerate}
\item[(1)] $\sigma(L^{op})=\upsilon(L^{op})$, i.e., $L$ is lean;
\item[(2)] $L^{op}$ is lean.
\end{enumerate}
\end{thm}
\noindent{\bf Proof.} 
$(1) \Rightarrow (2)$. Given any closed subset $F$ of $(L,\sigma(L))$, it is a compact saturated subset of $(L^{op},\upsilon(L^{op}))$ by Lemma \ref{closed compact op}. Thus, $F$ is compact saturated in $(L^{op},\sigma(L^{op}))$, i.e., $L^{op}$ is lean.

$(2) \Rightarrow (1)$. Assume that $L^{op}$ is lean. For any space $X$, denote $Q(X)$ the set of nonempty compact saturated subset of $X$ endowed with the reverse inclusion order. Given any closed subset $F$ of $(L,\sigma(L))$, it is compact saturated in $(L^{op},\sigma(L^{op}))$. By Lemma \ref{closed compact op}, $Q(\Sigma L^{op}) = Q((L^{op},\upsilon(L^{op})))$. Since $L$ is a continuous lattice, $\upsilon(L^{op}) = \omega(L)$ is continuous by Lemma \ref{hyper op}. Thus, $(L^{op},\upsilon(L^{op}))$ is core compact. Since $L^{op}$ is a complete lattice, $(L^{op},\upsilon(L^{op}))$  is locally compact and sober by Lemma \ref{up sober} and Lemma \ref{core sober}. 

We claim that $\Sigma L^{op}$ is core compact. Define $\eta:  L^{op} \to Q(\Sigma L^{op})$, $\eta(a) =\ \ua a$ and  $\Box:  \sigma(L^{op}) \to  \sigma (Q(\Sigma L^{op}))$, $\Box (U) = \{K \in Q(\Sigma L^{op}) : K \subseteq U \}$.  It is easy to see that $\eta$ is Scott continuous. $\eta^{-1}$ preserves arbitrary sups. $L^{op}$ is a complete lattice, so $\Sigma L^{op}$ is well-filtered. Then $\Box$ is well defined and Scott continuous.
$\eta^{-1} \circ \Box\ (U) =  \eta^{-1} (\{K \in Q(\Sigma L^{op}) : K \subseteq U \}) = U$. Thus $\eta^{-1} \circ \Box = id_{\sigma(L^{op})}$,  $\sigma( L^{op})$ is a Scott retract of $\sigma(Q(\Sigma L^{op}))$. Then, $\sigma( L^{op})$ is continuous, i.e., $\Sigma L^{op}$ is core compact. So it is  locally compact and sober by Lemma \ref{core sober}. 

By Hofmann-Mislove Theorem \cite[Theorem \uppercase\expandafter{\romannumeral2}-2.14]{CON03}, 
${\rm OFlit}(Q(X))$  is isomorphic to  $\mathcal{O}(X)$ for any locally compact sober space $X$ under the maps $g:{\rm OFilt}(Q(X)) \to \mathcal{O}(X)$, $g(\mathcal{F}) = \cup \mathcal{F}$ and $f:\mathcal{O}(X) \to {\rm OFilt}(Q(X))$, $f(U) = \Box U$. Then $(L^{op},\sigma( L^{op}))$ is equal to $(L^{op}, \upsilon(L^{op}))$. $\Box$

\vskip 3mm


\end{document}